\documentclass[reqno]{amsart}
\usepackage[english]{babel}
\usepackage{geometry,amssymb,enumerate,bbm,xcolor,amsthm}
\usepackage{cite}
\usepackage[T1]{fontenc}
\usepackage[utf8]{inputenc}
\usepackage{authblk}
\usepackage{lipsum}% http://ctan.org/pkg/lipsum
\makeatletter
\g@addto@macro{\endabstract}{\@setabstract}
\newcommand{\authorfootnotes}{\renewcommand\thefootnote{\@fnsymbol\c@footnote}}%
\makeatother

\catcode`\<=\active \def<{
	\fontencoding{T1}\selectfont\symbol{60}\fontencoding{\encodingdefault}}
\catcode`\>=\active \def>{
	\fontencoding{T1}\selectfont\symbol{62}\fontencoding{\encodingdefault}}

\newcommand{\mathd}{\mathrm{d}}

\newtheorem{corollary}{Corollary}
\newtheorem{definition}{Definition}
\newtheorem{lemma}{Lemma}
\newtheorem{proposition}{Proposition}
\newtheorem{remark}{Remark}%}
\newtheorem{theorem}{Theorem}
\begin{document}

\begin{center}
	\LARGE 
	New Regularity Criteria for the Navier-Stokes Equations in Terms of Pressure \par \bigskip
	
	\normalsize
	\authorfootnotes
	Benjamin Pineau\footnote{Corresponding author; Email:  bpineau@berkeley.edu}\textsuperscript{1}, Xinwei Yu\footnote{Email: xinwei2@ualberta.ca}\textsuperscript{2},
	\par \bigskip
	
	\textsuperscript{1}Department of Mathematics, University of California, Berkeley \par
	\textsuperscript{2}Department of Mathematics, University of Alberta \par \bigskip

	\today
\end{center}
\begin{abstract}
In this paper, we generalize the main results of {\cite{berselli2002regularity}} and {\cite{struwe2007serrin}} to Lorentz spaces, using a simple procedure. The main results are the following. Let $n\geq 3$ and let $u$ be a Leray-Hopf solution to the $n$-dimensional Navier-Stokes equations with viscosity $\nu$ and divergence free initial condition $u_0\in L^2(\mathbb{R}^n)\cap L^{k}(\mathbb{R}^n)$ (where $k=k(s)$ is sufficiently large). Then there exists a constant $c>0$ such that if
\begin{equation}
\|p\|_{L^{r,\infty}(0,\infty;L^{s,\infty}(\mathbb{R}^n))}<c\hspace{10mm}\frac{n}{s}+\frac{2}{r}\leq 2,\hspace{5mm}s>\frac{n}{2}
\end{equation}
or
\begin{equation}
\|\nabla p\|_{L^{r,\infty}(0,\infty;L^{s,\infty}(\mathbb{R}^n))}<c\hspace{10mm}\frac{n}{s}+\frac{2}{r}\leq 3,\hspace{5mm}s>\frac{n}{3}
\end{equation}
then $u$ is smooth on $(0, \infty) \times
\mathbb{R}^n$.{\hspace*{\fill}}{\medskip}
\\
Partial results in the case $n=3$ were obtained in \cite{suzuki2012regularity}, \cite{suzuki2012remark} and {\cite{1909.09960}}. Our results present a unified proof which works for all dimensions $n\geq 3$ and the full range or admissible pairs, $(s,r)$.
\end{abstract}

\section{Introduction}
The $n$-dimensional Navier-Stokes equations in $\mathbb{R}^n$ are given by
\begin{equation}
\partial_tu+(u\cdot\nabla)u=-\nabla p+\nu\Delta u, \label{eqn:1}
\end{equation}
\begin{equation}
\nabla\cdot u=0, \label{eqn:2}
\end{equation}
\begin{equation}
u(x,0)=u_0(x). \label{eqn:3}
\end{equation}
Herein $n\geq 3$, $x\in\mathbb{R}^n$, $t\in[0,T)$ for some $T>0$, $u:\mathbb{R}^n\mapsto\mathbb{R}^n$ is the velocity, $p:\mathbb{R}^n\mapsto\mathbb{R}$ is the pressure and $\nu>0$ is the viscosity. 

Denote by $H(\mathbb{R}^n)$, the $L^2$ closure of the space of divergence free, smooth, compactly supported functions on $\mathbb{R}^n$. It is a classical result of Leray {\cite{leray1934mouvement}} that for  $u_0\in H(\mathbb{R}^n)$, there exists a solution $u$ to (\ref{eqn:1}--\ref{eqn:3}) in the sense of distributions (commonly known as a Leray-Hopf solution), such that
\begin{equation}
u\in L^{\infty}(0,T;L^2(\mathbb{R}^n))\cap L^2(0,T;H^1(\mathbb{R}^n))
\end{equation}
Leray-Hopf solutions also satisfy the energy inequality,
\begin{equation}
\|u(t)\|_{L^2(\mathbb{R}^n)}^2+2\nu\int_{0}^{t}\|\nabla u(s)\|_{L^2(\mathbb{R}^n)}^2 ds\leq\|u_0\|_{L^2(\mathbb{R}^n)}^2
\end{equation}
for $t>0$. It has been known since \cite{leray1934mouvement} that smooth Leray-Hopf solutions are unique (in the sense that they are the only Leray-Hopf solutions satisfying their initial condition). In spite of much effort, the questions of global regularity and uniqueness of Leray-Hopf solutions remain open. Some partial progress  has been obtained in this direction and is discussed extensively in for instance, \cite{lemari2002recent,lemarie2018navier,robinson2016three,sohr2012navier}.

Although the question of global regularity remains open in general, one can place certain "integrability" assumptions on a Leray-Hopf solution $u$, to establish global regularity. For instance, it is a classical result of {\cite{prodi1959teorema,serrin1962interior}} that regularity of $u$ is guaranteed beyond some time $T>0$ as long as
\begin{equation}
\|u\|_{L^{r}(0,T;L^s(\mathbb{R}^n))}=\left(\int_{0}^{T}\|u\|_{L^s(\mathbb{R}^n)}^r dt\right)^{\frac{1}{r}}<\infty \label{eqn:4}
\end{equation}
where $s>n$ and $r$ satisfy
\begin{equation}
\frac{2}{r}+\frac{n}{s}\leq 1
\end{equation}
The above condition on $r$ and $s$ is significant because (\ref{eqn:4}) is invariant under the natural scaling of  (\ref{eqn:1}--\ref{eqn:3}), $u_\lambda(x,t):=\lambda u(\lambda x,\lambda^2 t)$. The case $s=n$, 
\begin{equation}
\|u\|_{L^{\infty}((0,T);L^n(\mathbb{R}^n))}<\infty,\label{eqn:ESS}
\end{equation} 
was established recently in {\cite{Escauriaza_2003}} (for $n=3$) and its later generalization \cite{DongDu2009}. 

The pioneering condition (\ref{eqn:4}) has inspired a large class of Prodi-Serrin type criteria. Instead of assumptions on $u$,  integrability conditions are placed on $\nabla u$, $p$, $\nabla p$, $\omega$, etc. Some of these conditions have also been further generalized using weaker norms. See for example, {\cite{da1995new}} ,{\cite{neustupa2002regularity}},{\cite{kozono2000bilinear}},{\cite{sohr2001regularity}}, {\cite{berselli2002regularity}}, {\cite{chan2007log}}, {\cite{qionglei2007regularity}}, {\cite{cao2008regularity}}, {\cite{fan2008regularity}}, {\cite{bersell2009some}}, {\cite{nunez2010regularity}},{\cite{vasseur2009regularity}}, {\cite{bjorland2011weak}}, {\cite{cao2011global}}, {\cite{tran2017regularity}}, {\cite{tran2017note}}, {\cite{han2019sharp}}, {\cite{pineau2018new}}, {\cite{chemin2017critical}}, {\cite{neustupa2018contribution}} and {\cite{pineau2020prodi}} for several results of this form.
\\
\\
One such generalization is due to Berselli and Galdi {\cite{berselli2002regularity}}. They proved the following.
\begin{theorem}
	Let $u$ be a
	Leray-Hopf solution to the Navier-Stokes equation with an initial value $u_0
	\in H(\mathbb{R}^n)\cap L^n (\mathbb{R}^n)$. If $p \in L^r (0, T ; L^s
	(\mathbb{R}^n))$ with
	\begin{equation}
	\frac{2}{r} + \frac{n}{s} \leq 2, \qquad s > \frac{n}{2}
	\end{equation}
	then $u$ is smooth on $(0, T] \times
	\mathbb{R}^n$, and can be extended beyond $T$.
\end{theorem}
An analogous theorem was also eventually proven for the pressure gradient. Partial results were obtained in {\cite{berselli2002regularity}}, {\cite{Zhou2005}} and eventually extended to the whole range of admissible pairs $(r,s)$ in {\cite{struwe2007serrin}}.
\\
The theorem from {\cite{struwe2007serrin}} is as follows.
\begin{theorem}
Let $u$ be a
Leray-Hopf solution to the Navier-Stokes equation with an initial value $u_0
	\in H(\mathbb{R}^n)\cap L^k (\mathbb{R}^n)$. If $\nabla p \in L^r (0, T ; L^s
(\mathbb{R}^n))$ with
\begin{equation}
\frac{2}{r} + \frac{n}{s} \leq 3, \qquad s > \frac{n}{3}
\end{equation}
then $u$ is smooth on $(0, T]\times
\mathbb{R}^n$ and can be extended beyond $T$.{\hspace*{\fill}}{\medskip}
\end{theorem}
\begin{remark}
Note that in the assumption on the initial condition, we generally have $k=k(s)$. This is at least the case for $n>4$ in light of the approach of {\cite{struwe2007serrin}}. For $n=3,4$, one may take $k=4$ in light of the method of {\cite{Zhou2005}}.
\end{remark}
	In this paper, we aim to extend the above theorems to Lorentz spaces. In particular, we prove the following results and corollaries. Herein, we also assume $u$ is a Leray-Hopf solution to (\ref{eqn:1}--\ref{eqn:3}) with pressure $p$ and also that $u_0\in H(\mathbb{R}^n)\cap L^{k}(\mathbb{R}^n)$
	where $k=k(s)$ is chosen to guarantee sufficient spatial integrability in the following proofs. 
	
	\begin{theorem}\label{thm:1}
	There exists a constant $c>0$ such that if
	\begin{equation}
	\|p\|_{L^{r,\infty}(0,T;L^{s,\infty}(\mathbb{R}^n))}\leq c\hspace{10mm}\frac{n}{s}+\frac{2}{r}\leq 2,\hspace{5mm}s>\frac{n}{2}\label{eqn:5}
	\end{equation}
	then $u$ is smooth on $(0, T] \times
	\mathbb{R}^n$ and can be extended beyond $T$.{\hspace*{\fill}}{\medskip}
	\end{theorem}
\begin{theorem}{\label{thm:gradp}}
	There exists a constant $c>0$ such that if
	\begin{equation}
	\|\nabla p\|_{L^{r,\infty}(0,T;L^{s,\infty}(\mathbb{R}^n))}\leq c\hspace{10mm}\frac{n}{s}+\frac{2}{r}\leq 3,\hspace{5mm}s>\frac{n}{2}
	\end{equation}
	then $u$ is smooth on $(0, T] \times
	\mathbb{R}^n$ and can be extended beyond $T$.{\hspace*{\fill}}{\medskip}
\end{theorem}
Partial results on regularity criteria for weak-Lebesgue spaces for the pressure and pressure gradient have been obtained in  the case $n=3$ in {\cite{suzuki2012remark}}, {\cite{suzuki2012regularity}} and {\cite{1909.09960}}. Our approach presents a unified method of handling all pairs of appropriate $(s,r)$ in any dimension $n\geq 3$
\\
\\
We can also strengthen the assumptions of (\ref{eqn:5}) to obtain the following.
\begin{corollary}\label{thm:3}
Suppose we have
\begin{equation}
\|p\|_{L^{r,r'}(0,T;L^{s,\infty}(\mathbb{R}^n))}<\infty
\end{equation}
for $s>n/2$ and $r'\in (0,\infty)$ satisfying $\frac{n}{s}+\frac{2}{r}\leq 2$, then $u$ is smooth on $(0, T] \times
\mathbb{R}^n$ and can be extended beyond $T$.{\hspace*{\fill}}{\medskip}
\end{corollary}
Similarly, we have
\begin{corollary}
Suppose we have
\begin{equation}
\|\nabla p\|_{L^{r,r'}(0,T;L^{s,\infty}(\mathbb{R}^n))}<\infty
\end{equation}
for $s>n/3$ and $r'\in (0,\infty)$ satisfying $\frac{n}{s}+\frac{2}{r}\leq 3$, then $u$ is smooth on $(0, T] \times
\mathbb{R}^n$ and can be extended beyond $T$.{\hspace*{\fill}}{\medskip}
\end{corollary}
\begin{remark}
	In light of Proposition \ref{prop:1}, these two corollaries generalize Theorems 1 and 2.
\end{remark}

\begin{remark}
    Upon the completion of this paper, a new preprint \cite{1909.09960} was posted on arXiv, where the case $n=3$ of Theorems \ref{thm:1} and \ref{thm:gradp} are proved using similar ideas. We note that the proofs of the general case requires new ideas and techniques, especially for Theorem \ref{thm:gradp}. 
\end{remark}
	
\section{Preliminaries}
For convenience, we provide the definition and some basic properties of Lorentz spaces. See {\cite{grafakos2008classical}} for a more detailed discussion.
\begin{definition} If $(X,\mu)$ is some measure space and $0<p<\infty$, $0<q\leq\infty$, then the Lorentz space $L^{p,q}(X,\mu)$ is the space of all measurable functions for which the quasinorm,
\begin{equation}
\begin{split}
\|f\|_{L^{p,q(X)}}:=\begin{cases} p^{\frac{1}{q}}\left(\int_{0}^{\infty}(d_f(s)^{\frac{1}{p}}s)^q\frac{ds}{s}\right)^{\frac{1}{q}},\hspace{4mm}&q<\infty
\\
\sup_{s>0}sd_f(s)^{\frac{1}{p}},&q=\infty
\end{cases}
\end{split}
\end{equation}
is finite. Here $d_f$ is defined by
\begin{equation}
\begin{split}
d_f(s):=\mu(\{x\in X:|f(x)|>s\})
\end{split}
\end{equation}
\end{definition}
We also recall some basic properties of Lorentz spaces.
\begin{proposition}[Basic Properties of Lorentz Spaces]~\\
\label{prop:1}
\begin{enumerate}
    \item 
    $L^p(X,\mu)=L^{p,p}(X,\mu)$
    \item $\||f|^r\|_{L^{p,q}(X,\mu)}=\|f\|_{L^{pr,qr}(X,\mu)}^r$ for $0<p,r<\infty$ and $0<q\leq\infty$
    \item $\|f\|_{L^{p,q}(X,\mu)}\leq C\|f\|_{L^{p,r}(X,\mu)}$ for $0<p<\infty$ and $0<q_2\leq q_1\leq\infty$
\end{enumerate}
\end{proposition}
Lorentz spaces also enjoy a variant of Young's inequality.
\begin{proposition}[Young's Inequality for Lorentz Spaces {\cite{lemari2002recent}}]
\label{prop:2}
Let $1<p<\infty$, $1\leq q\leq\infty$ and $\frac{1}{p'}+\frac{1}{p}=1$, $\frac{1}{q'}+\frac{1}{q}=1$. Suppose as well that $1<p_1<p'$ and $q'\leq q\leq\infty$. If $\frac{1}{p_2}+1=\frac{1}{p}+\frac{1}{p_1}$ and $\frac{1}{q_2}=\frac{1}{q}+\frac{1}{q_1}$, then the convolution operator,
\end{proposition}
\begin{equation}
*:L^{p,q}(\mathbb{R}^n)\times L^{p_1,q_1}(\mathbb{R}^n)\mapsto L^{p_2,q_2}(\mathbb{R}^n)
\end{equation}
is a bounded bilinear operator.
\\
\\
A consequence of the above proposition is the following stronger variant of the Sobolev inequality for $n>2$,
\begin{equation}
\|u\|_{L^{\frac{2n}{n-2},2}(\mathbb{R}^n)}\leq C\|\nabla u\|_{L^2(\mathbb{R}^n)}
\end{equation}
Another useful consequence is the following Sobolev-type inequality for weak $L^p$ spaces.That is, for $\frac{1}{p^{*}}=\frac{1}{p}-\frac{1}{n}$ and $1<p<n$, we have
\begin{equation}
\|u\|_{p^{*},\infty}\leq C\|\nabla u\|_{p,\infty}
\end{equation}
There is also a variant of H{\"o}lder's inequality for Lorentz spaces, due to O'Neil.
\begin{proposition}[H{\"o}lder's Inequality for Lorentz Spaces {\cite{o1963convolution}}]
\label{prop:3}
Let $0<p_1,p_2,p<\infty$ and $0<q_1,q_2\leq\infty$ satisfy $\frac{1}{p}=\frac{1}{p_1}+\frac{1}{p_2}$ and $\frac{1}{q}=\frac{1}{q_1}+\frac{1}{q_2}$. Then for $f\in L^{p_1,q_1}$ and $g\in L^{p_2,q_2}$ we have,
\begin{equation}
\|fg\|_{L^{p,q}(X)}\leq C\|f\|_{L^{p_1,q_1}(X)}\|g\|_{L^{p_2,q_2}(X)}.
\end{equation}
The constant $C$ depends only on $p_1,q_1,p_2,q_2$.
\end{proposition}
The following sharper version of a standard interpolation inequality will be crucial in our proof of Theorem \ref{thm:gradp} for dimensions $n\geqslant 4$. See Proposition 1.14 of {\cite{grafakos2008classical}} for details.
\begin{proposition}
Let $0<p<q\leq\infty$ and let $f\in L^{p,\infty}(\mathbb{R}^n)\cap L^{q,\infty}(\mathbb{R}^n)$. Then $f\in L^{r}(\mathbb{R}^n)$ for $p<r<q$ and we have
	\begin{equation}
	\|f\|_{L^r}\leq C\|f\|_{L^{p,\infty}}^{\frac{\frac{1}{r}-\frac{1}{q}}{\frac{1}{p}-\frac{1}{q}}}\|f\|_{L^{q,\infty}}^{\frac{\frac{1}{p}-\frac{1}{r}}{\frac{1}{p}-\frac{1}{q}}}
	\end{equation}
	The constant $C$ depends only on $p,r,q$.
\end{proposition}

We recall a Lorentz space variant (see {\cite{lemari2002recent}}) of a standard  Calder\'on-Zygmund type inequality that allows us to estimate the pressure in terms of the velocity.

\begin{proposition}
For $1<p<\infty$, we have
\begin{equation}
\|p\|_{L^{p,q}(\mathbb{R}^n)}\leq C\|u\|_{L^{{2p,2q}}(\mathbb{R}^n)}^2
\end{equation}
\end{proposition}

Finally, the following ``nonlinear Gronwall-type Inequality'', established in {\cite{pineau2020prodi}}, is crucial in the proofs of our theorems. 

\begin{lemma}
\label{lem:1}
Let $T>0$ and $\varphi\in L_{loc}^{\infty}([0,T))$ be non-negative. Assume further that
\begin{equation}
\begin{split}
\varphi (t)\leq C_0+C_1\int_{0}^{t}\mu(s)\varphi(s)ds+\kappa\int_{0}^{t}\lambda(s)^{1-\epsilon}\varphi(s)^{1+A(\epsilon)}ds\hspace{5mm}\forall 0<\epsilon<\epsilon_0.
\end{split}
\end{equation}
Where $\kappa,\epsilon_0 >0$ are constants, $\mu\in L^1(0,T)$ and $A(\epsilon)>0$ satisfies $\lim_{\epsilon\to 0}\frac{A(\epsilon)}{\epsilon}=c_0>0$. Then $\varphi$ is bounded on $[0,T]$ if $\|\lambda\|_{L^{1,\infty}(0,T)}<c_0^{-1}\kappa^{-1}$
\end{lemma}

\section{Proofs of Main Results}

\subsection{Proof of Theorem \ref{thm:1}}
It clearly suffices to prove the theorem for $s,r$ satisfying $\frac{n}{s}+\frac{2}{r}=2$. Herein, we assume $s,r$ have this property.
\\
\\
There are two cases to consider.
\begin{itemize}
\item Case 1. $\frac{n}{2}<s\leq\frac{n+2}{2}$.
In this case we have
\begin{equation}
\begin{split}
2s\leq\frac{ns}{s-1}<\frac{2ns}{n-2}
\end{split}
\end{equation}
Taking the $L^2$ inner product of $u|u|^{2s-2}$ and (\ref{eqn:1}) and integrating by parts yields
\begin{equation}
\begin{split}
\frac{d}{dt}\|u\|_{2s}^{2s}+\|\nabla |u|^s\|_2^2 & \leq\int_{\mathbb{R}^n}|p|^{\frac{1}{2}}|p|^{\frac{1}{2}}|u|^{s-1}|\nabla |u|^s|
\\
&\leq \||p|^{\frac{1}{2}}\|_{L^{2s,\infty}}\||p|^{\frac{1}{2}}\|_{L^{\frac{ns}{s-1},2s}}\||u|^{s-1}\|_{L^{\frac{2ns}{(s-1)(n-2)},\frac{2s}{s-1}}}\|\nabla |u|^s\|_{L^2}
\\
&= \||p|^{\frac{1}{2}}\|_{L^{2s,\infty}}\||p|^{\frac{1}{2}}\|_{L^{\frac{ns}{s-1},2s}}\||u|^s\|_{L^{\frac{2n}{n-2},2}}^{1-\frac{1}{s}}\|\nabla |u|^s\|_{L^2}
\\
&\leq C\||p|^{\frac{1}{2}}\|_{L^{2s,\infty}}\|u\|_{L^{\frac{ns}{s-1},2s}}\|\nabla |u|^s\|_{L^2}^{2-\frac{1}{s}}
\\
&\leq C\||p|^{\frac{1}{2}}\|_{L^{2s,\infty}}\|u\|_{L^{2s,2s}}^{s-\frac{n}{2}}\|u\|_{L^{\frac{2ns}{n-2},2s}}^{\frac{n}{2}+1-s}\|\nabla |u|^s\|_{L^2}^{2-\frac{1}{s}}
\\
&\leq C\|p\|_{L^{s,\infty}}^{\frac{1}{2}}\|u\|_{L^{2s,2s}}^{s-\frac{n}{2}}\|\nabla |u|^s\|_{L^2}^{\frac{2s+n}{2s}}
\\
&\leq C\|p\|_{L^{s,\infty}}^{\frac{2s}{2s-n}}\|u\|_{L^{2s,2s}}^{2s}+\frac{1}{2}\|\nabla |u|^s\|_{L^2}^{2}
\end{split}
\end{equation}
It follows that
\begin{equation}
\begin{split}
\frac{d}{dt}\|u\|_{2s}^{2s}+\|\nabla |u|^s\|_2^2 &\leq  C\|p\|_{L^{s,\infty}}^{\frac{2s}{2s-n}}\|u\|_{L^{2s,2s}}^{2s}
\\
&\leq C\|p\|_{L^{s,\infty}}^{\frac{2s}{2s-n}(1-\epsilon)}\|u\|_{L^{2s,2s}}^{2s(1+\frac{2\epsilon}{2s-n})}
\end{split} \label{eqn:6}
\end{equation}

Integrating (\ref{eqn:6}) from $0$ to $t<T$ and applying Lemma \ref{lem:1} shows that $\|u\|_{2s}$ is bounded on $[0,T]$ provided (\ref{eqn:5}) is satisfied for some small enough constant, $c>0$. In particular, this means $u\in L^{r}(0,T;L^{2s}(\mathbb{R}^n))$ for $r$ satisfying
\begin{equation}
\begin{split}
\frac{n}{2s}+\frac{2}{r}=1
\end{split}
\end{equation}
Regularity of $u$ then follows from the classical Prodi-Serrin condition for the velocity $u$.
\\
\\
\item Case 2: $s>\frac{n+2}{2}$
\\
\\
In this case we have
\begin{equation}
\begin{split}
2s<\frac{2s(s-1)}{s-2}<\frac{2ns}{n-2}
\end{split}
\end{equation}
As in the first case, taking the $L^2$ inner product of $u|u|^{2s-2}$ and  (\ref{eqn:1}) and integrating by parts yields
\begin{equation}
\begin{split}
\frac{d}{dt}\|u\|_{L^{2s}(\mathbb{R}^n)}^{2s}+\|\nabla |u|^s\|_{L^2}^2 &\leq\int_{\mathbb{R}^n}|p||u|^{s-1}|\nabla |u|^s|
\\
&\leq C\||p|^{\frac{1}{2}}\|_{L^{2s,\infty}}\||p|^{\frac{1}{2}}\|_{L^{2s,2s}}\||u|^{s-1}\|_{L^{\frac{2s}{s-2},\frac{2s}{s-1}}}\|\nabla |u|^s\|_{L^2}
\\
&= C\|p\|_{L^{s,\infty}}^{\frac{1}{2}}\|p\|_{L^{s}}^{\frac{1}{2}}\|u\|_{L^{\frac{2s(s-1)}{s-2},2s}}^{s-1}\|\nabla |u|^s\|_{L^2}
\\
&\leq C\|p\|_{L^{s,\infty}}^{\frac{1}{2}}\|p\|_{L^{s}}^{\frac{1}{2}}\|u\|_{L^{2s,2s}}^{s-1-\frac{n}{2}}\|u\|_{L^{\frac{2ns}{n-2},2s}}^{\frac{n}{2}}\|\nabla |u|^s\|_{L^2}
\end{split}
\end{equation}
Then, by Young's inequality we obtain,
\begin{equation}
\begin{split}
\frac{d}{dt}\|u\|_{L^{2s}(\mathbb{R}^n)}^{2s}+\|\nabla |u|^s\|_{L^2(\mathbb{R}^n)}^2 &\leq C\|p\|_{L^{s,\infty}}^{\frac{2s}{2s-n}}\|p\|_{L^{s}}^{\frac{2s}{2s-n}}\|u\|_{L^{2s,2s}}^{\frac{2s(2s-n-2)}{2s-n}}
\\
&\leq C\|p\|_{L^{s,\infty}}^{\frac{2s}{2s-n}}\|u\|_{L^{2s}}^{\frac{4s}{2s-n}}\|u\|_{L^{2s,2s}}^{\frac{2s(2s-n-2)}{2s-n}}
\\
&=C\|p\|_{L^{s,\infty}}^{\frac{2s}{2s-n}}\|u\|_{L^{2s}}^{2s}
\\
&\leq C\|p\|_{L^{s,\infty}}^{\frac{2s}{2s-n}(1-\epsilon)}\|u\|_{L^{2s}}^{2s(1+\frac{2\epsilon}{2s-n})}
\end{split}
\end{equation}
Regularity of $u$ now follows from the same approach used in case 1.
\end{itemize}

\subsection{Proof of Corollary 1}
We proceed in a similar way to {\cite{sohr2012navier}}. From the proof of Theorem \ref{thm:1}, we have
\begin{equation}
\label{eqn:8}
\begin{split}
\frac{d}{dt}\|u\|_{L^{2s}(\mathbb{R}^n)}^{2s} \leq C\|p\|_{L^{s,\infty}}^{\frac{2s}{2s-n}(1-\epsilon)}\|u\|_{L^{2s}}^{2s(1+\frac{2\epsilon}{2s-n})}
\end{split}
\end{equation}
Now, we observe that for $0<t_0<T$ and $f\in L^{p,r}(t_0,T)$ there holds,
\begin{equation}
\begin{split}
\| f \|_{L^{r, r'} (t_0, T)} & = r^{1 / r'}  \left( \int_0^{\infty}
(d_f (s)^{1 / r} s)^{r'}  \frac{\mathd s}{s} \right)^{1 / r'} \nonumber
\\
& \leq r^{1 / r'}  \left[ \left( \int_0^R (d_f (s)^{1 / r} s)^{r'}
\frac{\mathd s}{s} \right)^{1 / r'} + \left( \int_R^{\infty} (d_f (s)^{1 / r'}
s)^{r'}  \frac{\mathd s}{s} \right)^{1 / r'} \right] \nonumber
\\
& \leq r^{1 / r'}  \left[ (T-t_0)^{1 / r}\left(\int_0^R s^{r' - 1}
\mathd s \right)^{1 / r'} + \left( \int_R^{\infty} (d_f (s)^{1 / r'} s)^{r'} 
\frac{\mathd s}{s} \right)^{1 / r'} \right] . 
\end{split}
\end{equation}
We see in light of Proposition \ref{prop:1} and by taking $R$ sufficiently large and then $t_0$ close enough to $T$, we can make $\|p\|_{L^{r,\infty}(t_0,T;L^s(\mathbb{R}^n))}$ arbitrarily small. Thus, the claim follows from Theorem \ref{thm:1}.

\subsection{Proof of Theorem \ref{thm:gradp}}
It clearly suffices to prove the theorem for $s,r$ satisfying $\frac{n}{s}+\frac{2}{r}=3$. Herein, we assume $s,r$ have this property. Again, we have two cases. 
\\
\begin{itemize}
\item Case 1: $\frac{n}{3}<s<n$
\\
\\
Recall that from the proof of theorem 3, we have for any $\theta>\frac{n}{2}$,
\begin{equation}
\begin{split}
\frac{d}{dt}\|u\|_{2\theta}^{2\theta}+\|\nabla |u|^\theta\|_2^2 &\leq  C\|p\|_{L^{\theta,\infty}}^{\frac{2 \theta}{2\theta-n}}\|u\|_{L^{2\theta,2\theta}}^{2\theta}
\\
&\leq C\|p\|_{L^{\theta,\infty}}^{\frac{2\theta}{2\theta-n}(1-\epsilon)}\|u\|_{L^{2\theta,2\theta}}^{2\theta(1+\frac{2\epsilon}{2\theta-n})}
\end{split} 
\end{equation}
Taking
$\theta=\frac{ns}{n-s}>\frac{n}{2}$ and using the Sobolev-type inequality,
\begin{equation}
\|p\|_{\theta,\infty}\leq C\|\nabla p\|_{s,\infty}
\end{equation}
yields the desired conclusion.
\\
\item Case 2: $s\geq n$
\\
\\
This case is more technical, but the main idea is similar to the proof of Theorem 3. Indeed, for some $\theta>0$ sufficiently large (i.e. $\theta>s$), we can multiply (\ref{eqn:1}) by $u|u|^{(\theta-2)}$ and integrate by parts to obtain
\begin{equation}
\begin{split}
\frac{d}{dt}\|u\|_{\theta}^{\theta}+\|\nabla |u|^{\frac{\theta}{2}}\|_2^2 \leq C\left|\int_{\mathbb{R}^n}u\cdot\nabla p|u|^{\theta-2}dx\right|:=CI
\end{split}
\end{equation}
Inspired by the approach of {\cite{struwe2007serrin}}, we estimate $I$ in two different ways. Indeed, we can bound $I$ by
\begin{equation}
\begin{split}
I\leq I_1:=\int_{\mathbb{R}^n}|\nabla p||u|^{\theta-1}\ \mathd x
\end{split}
\end{equation}
and
\begin{equation}
\begin{split}
I\leq I_2:=\int_{\mathbb{R}^n}|p||\nabla |u|||u|^{\theta-2}\ \mathd x
\end{split}
\end{equation}
\textbf{Estimate for $I_1$:}
We have for some $\delta=\delta(\theta,s)>0$ to be chosen later, and Proposition 4,
\begin{equation}
\begin{split}
I_1=\||\nabla p||u|^{\theta-1}\|_{1}&=\||\nabla p|^{\frac{1}{2}}|u|^{\frac{\theta-1}{2}}\|_{2}^2
\\
&\leq C\||\nabla p|^{\frac{1}{2}}|u|^{\frac{\theta-1}{2}}\|_{2-\delta,\infty}\||\nabla p|^{\frac{1}{2}}|u|^{\frac{\theta-1}{2}}\|_{\frac{2-\delta}{1-\delta},\infty}
\end{split}
\end{equation}
Taking $\delta=\frac{2(\theta-s)}{\theta(s+1)-s}$ and apply Holder's inequality for Lorentz spaces gives
\begin{equation}
\begin{split}
I_1&\leq C\||\nabla p|^{\frac{1}{2}}\|_{2s,\infty}\||\nabla p|^{\frac{1}{2}}\|_{2s,\infty}\|u\|_{\theta,\infty}^{\frac{\theta-1}{2}}\|u\|_{\frac{\theta (\theta-1)s}{\theta s-2\theta+s},\infty}^{\frac{\theta-1}{2}}
\\
&=C\|\nabla p\|_{s,\infty}\|u\|_{\theta,\infty}^{\frac{\theta-1}{2}}\|u\|_{\frac{\theta (\theta-1)s}{\theta s-2\theta+s},\infty}^{\frac{\theta-1}{2}}
\\
&\leq C\|\nabla p\|_{s,\infty}\|u\|_{\theta,\infty}^{\frac{\theta-1}{2}}\|u\|_{\theta,\infty}^{\frac{\theta s+ns-\theta n-s}{2s}}\|u\|_{\frac{n}{n-2}\theta}^{\frac{n(\theta-s)}{2s}}
\\
&\leq C\|\nabla p\|_{s,\infty}\|u\|_{\theta,\infty}^{\frac{\theta-1}{2}+\frac{\theta s+ns-\theta n-s}{2s}}\|\nabla |u|^{\frac{\theta}{2}}\|_2^{\frac{n(\theta-s)}{\theta s}}
\end{split}
\end{equation}
\textbf{Estimate for $I_2$}: 
\begin{equation}
\begin{split}
I_2 =\int_{\mathbb{R}^n}|p||\nabla|u|||u|^{\theta-2}&\leq \left(\int_{\mathbb{R}^n}|p|^2|u|^{\theta-2}\right)^{\frac{1}{2}}\|\nabla |u|^{\frac{\theta}{2}}\|_2
\\
&\leq \|p\|_{\frac{\theta+2}{2}}\|u\|_{\theta+2}^{\frac{\theta-2}{2}}\|\nabla |u|^{\frac{\theta}{2}}\|_2
\\
&\leq C\|u\|_{\theta+2}^{\frac{\theta+2}{2}}\|\nabla |u|^{\frac{\theta}{2}}\|_2
\\
&\leq C\|u\|_{\theta}^{\frac{\theta-n+2}{2}}\|\nabla |u|^{\frac{\theta}{2}}\|_2^{\frac{n}{\theta}+1}
\end{split}
\end{equation}
Now, we estimate $I$ by combining the above two estimates. We have by Young's inequality, and after possibly relabeling $\epsilon$,
\begin{equation}
\begin{split}
I&\leq I_1^{\frac{1}{2}-\epsilon}I_2^{\frac{1}{2}+\epsilon}
\\
&\leq C\|\nabla p\|_{s,\infty}^{\frac{1}{2}-\epsilon}\|u\|_{\theta}^{(\frac{1}{2}-\epsilon)\frac{\theta-1}{2}+(\frac{1}{2}-\epsilon)\frac{\theta s+ns-\theta n-s}{2s}+(\frac{1}{2}+\epsilon)\frac{\theta-n+2}{2}}\|\nabla |u|^{\frac{\theta}{2}}\|_2^{(\frac{1}{2}-\epsilon)\frac{n(\theta-s)}{\theta s}+(\frac{n}{\theta}+1)(\frac{1}{2}+\epsilon)}
\\
&\leq C\|\nabla p\|_{s,\infty}^{\frac{2s}{3s-n}(1-c_1\epsilon)}\|u\|_{\theta}^{\theta+c_2\epsilon+O(\epsilon^2)}+\|\nabla |u|^{\frac{\theta}{2}}\|_2^{2}
\end{split}
\end{equation}
where $c_1$ and $c_2$ are positive constants. From this, we conclude
\begin{equation}
\begin{split}
\frac{d}{dt}\|u\|_{\theta}^{\theta}\leq C\|\nabla p\|_{s,\infty}^{\frac{2s}{3s-n}(1-c_1\epsilon)}\|u\|_{\theta}^{\theta+c_2\epsilon+O(\epsilon^2)}
\end{split}
\end{equation}
So, the claim follows by the nonlinear Gronwall lemma.
\end{itemize}
\begin{remark}
We remark that the above treatment is necessary in dealing with higher dimensions $n\geq 4$. The proof in {\cite{1909.09960}} takes $\theta=4$ which is crucial in their being able to take advantage of the relation $div(\nabla p + u\cdot\nabla u)=0$. As $\theta\geq n$ is necessary to conclude regularity, $\theta=4$ is sufficient for $n=3,4$, but not for larger $n$.
\end{remark}
\subsection{Proof of Corollary 2}
The argument is essentially identical to the proof of Corollary 1, so we omit the proof.

\bibliographystyle{plain}
\bibliography{Pressure_paper}
\end{document}